%
%
%
%
%
%

\documentstyle{amsppt}

\magnification =1200
\def \headd{ \headline={\sevenrm \hfil \jobname    
\hfil\number\month/ \number\day/\number\year \ \ \number\time}}
\def\cn {{ \Bbb C}^{N}}
\def \zzb {\bar Z}
\def \prt#1 { {\partial \over \partial {#1   }}  }
\def \prtf #1 #2{ {\partial #1 \over \partial {#2   }  }   }

\def \zo {Z_0}
\def \po {p_0}

\def \zb {\bar z}
\def \wb {\bar w}
\def \z {\zeta}
\def \t {\tau}
\def \a {\alpha}

\def \co {\Cal O}

\def \Z {\bar Z}

\def \no {{n+1}}
\def \nt {{n+2}}
\def \H {\bar H}
\def \h {\text}
\def\-{\overline}
\def \od  {\hbox {\rm ord}}

\def \gb {\bar g}
\def \fb {\bar f}

\def \zo {z^0}
\def \wo {w^0}

\def \zob {\zb^0}

\def \b {\beta}
\def \zzb {\bar Z}
\def \cn {\Bbb C^{N }} 
\def \spo {\Sigma_{\po}}
\def \Qb {{\bar Q}}
\def \ch {\chi}
\def \pob {\bar \po}
\def \zob {\bar \zo}
\def \de  {\Delta_\epsilon}
\def \no {{n+1}}
\def \nt {{n+2}}
\def \H {\bar H}
\def \h {\text}
\def \od  {\hbox {\rm ord}}


\document

\topmatter
\author  M.S. Baouendi, Xiaojun Huang, and Linda Preiss Rothschild
\endauthor
\title Regularity of CR mappings between algebraic hypersurfaces
\endtitle
\address  M.S.B.and L.P.R. :Department of
Mathematics-0112, University of California, \hfill\break San  
Diego, La Jolla,
CA 92093-0112,\newline X.H.: Department of Mathematics, University
of Chicago, Chicago, IL 60637
\endaddress
\email sbaouendi\@ucsd.edu, xhuang\@math.uchicago.edu,
lrothschild\@ucsd.edu \endemail
\thanks The first and third authors were partially supported by
National Science Foundation Grant DMS 9203973  \endthanks
\rightheadtext{Regularity of CR mappings}
\leftheadtext {M.S. Baouendi, X. Huang, and L. P.
Rothschild }
\endtopmatter 
\heading  \bf \S 0 Introduction
\endheading

Let $M$ and $M'$ be real analytic hypersurfaces in $\Bbb C^{ N}$
and $\Bbb C^{ N'}$ respectively and $H:M \to M'$ a sufficiently
smooth CR mapping. Under what conditions does $H$ extend
holomorphically to a neighborhood of $M$ in $\Bbb C^{N }$?  In
this paper we prove that if $M$ and $M'$ are algebraic
hypersurfaces in $\Bbb C^{ N}$, i.e. both defined by the
vanishing of real polynomials, then any sufficiently
smooth CR mapping with Jacobian not identically zero extends
holomorphically provided the hypersurfaces are holomorphically
nondegenerate (see definition below). Conversely, we prove
that holomorphic nondegeneracy is necessary for this property
of CR mappings to hold. For the case of unequal dimensions, we
also prove that if
$N'=N+1$, $M'$ is the sphere, and $M$ is an algebraic hypersurface
which does not contain any complex variety of positive codimension,
extendability holds for all CR mappings with certain minimal a priori
regularity. 

Our approach uses the work of Webster [W1], [W2], on
holomorphic mappings between algebraic hypersurfaces, and the
recent generalizations in [H1], [H2] and [BR6].  The question
of holomorphic extendability of CR mappings between real
analytic hypersurfaces has attracted considerable attention
since the work of Lewy [Lw] and Pincuk [P]. For more recent
work in the case $N = N'$, see Diederich-Webster
[DW], Jacobowitz,Treves, and the first author [BJT], Bell and
the first and third authors [BBR], [BR1], [BR2],
Diederich-Fornaess [DF], [BR3], and the references therein, as
well as the survey paper Forstneri\v c [Fo2]. We note here that 
the results for
$N\ge 3$ cited above require  nonvanishing conditions on
the normal component of the mapping and require the first
hypersurface to be  essentially finite.   
(See [BR3] for more general results for the case $N=2$ and Meylan
[Me1], [Me2]  for some extensions of this to higher dimensions.)  
In the algebraic case, studied in this paper,  we are able to omit
these assumptions. A recent example given by Ebenfelt [E] shows
that holomorphic extendibility may fail if the hypersurfaces are
not assumed to be algebraic. The authors know of
no other example of real analytic hypoellipticity which holds in
the algebraic category but not in the real analytic category.

 For the
case where
$N
\not= N'$ an important first result was given by Webster
[W2], who proved that any CR map of class
$C^3$ from a strongly pseudoconvex real analytic hypersuface in
$\Bbb C^{N }$ to the sphere in $\Bbb C^{N+1 }$ admits a
holomorphic extension  on a dense open subset. Generalizations
were later given by Faran [Fa1] [Fa2],
Cima-Suffridge [CS1] [CS2], Cima-Krantz-Suffridge [CKS], 
Forstneri\v c [Fo1], and [Hu2].  
Recently, the second author in [H1], [H2] proved
that any CR mapping of class $C^{N'-N+1}$
between two strictly pseudoconvex real analytic hypersurfaces in
$\Bbb C^{N }$ and $\Bbb C^{N' }$  ($N'\ge N>1$) respectively,  
is real analytic on a dense open subset of $M$, and is algebraic
if both
$M$ and $M'$ are algebraic.
In Theorem 5 below we prove that holomorphic extension holds
everywhere under weaker differentiability assumptions than those
given in [H2]. 

We now introduce some notation and definitions which are needed
to state precisely  our main results.   By a germ at $\po$ of a holomorphic vector field in
$\cn$, we shall mean a complex vector field  of the form  $\sum_1^N
a_j(Z)\prt {Z_j}
$, where the $a_j(Z)$ are germs at $\po$ of holomorphic functions.
 Let
$M$ be a real analytic hypersurface in
$\Bbb C^{N }$. For $\po \in M$ we say that  $M$ is {\it
holomorphically degenerate} at $\po$  if there exists a nonzero
germ of a holomorphic vector field tangent to
$M$ in a neighborhood of
$\po$ (see Stanton [Sta], [BR6]). We say that  $M$ is {\it
holomorphically nondegenerate}  if it is not holomorphically
degenerate at any $p_0$ in
$M$. Recall that by Theorem 1 of [BR6], a connected real
analytic hypersurface is holomorphically nondegenerate if and
only if there is a point $p_1$ at which it is not holomorphically
degenerate.  A CR function on $M$ is a function which is
annihilated by the tangential Cauchy-Riemann operators; a mapping
from $M$ into $\Bbb C^{ N'}$ is CR if its components are CR
functions.   

\proclaim {Theorem 1} Let $M$ and $M'$ be two algebraic
hypersurfaces in $\Bbb C^{ N}$ and assume that $M$ is connected
and holomorphically nondegenerate.  If $H$ is a smooth CR mapping
from $M$ to $M'$ with  {\rm Jac }$H \not\equiv 0$, where {\rm Jac
}$H$ is the Jacobian determinant of $H$, then $H$ extends
holomorphically in an open neighborhood of $M$ in $\Bbb C^{ N}$. 
\endproclaim 

The fact that $M$ and $M'$ are algebraic plays an important role.
 Indeed, as mentioned before, a recent example given by  Ebenfelt [E]
shows that the conclusion of Theorem 1 need not hold if $M$ is real
analytic, but not algebraic.  (See Example 2.10 below.)

Note that if the Jacobian of a CR map is 0, then $M'$ must contain
a complex variety (see [BR5]).  Therefore we obtain the following
corollary of Theorem 1.
  
\proclaim {Corollary 1} Let $M$ and $M'$ be two algebraic
hypersurfaces in $\Bbb C^{ N}$ and assume that $M$ is connected
and holomorphically nondegenerate.  If $H$ is a smooth CR mapping
from $M$ to $M'$ and if $M'$ contains no complex analytic variety
of positive dimension, then
$H$ extends holomorphically in an open neighborhood of $M$ in $\Bbb
C^{ N}$. 
\endproclaim

If $f$ is a function defined on $M$ we shall say that $f$ is
{\it algebraic} if there exist holomorphic polynomials $q_j(Z), j =
0,\ldots, k$, not all identically $0$, such that $q_k(Z)f(Z)^k +
\ldots + q_0(Z) \equiv 0$, for $Z \in M$.  Similarly, we say that
$f$ is {\it locally algebraic} if for any point $p$ on $M$ there is
a neighborhood of $p$ such that the restriction of $f$ to that
neighborhood is algebraic.  A mapping is algebraic (resp. locally
algebraic) if each of its components is.  

\proclaim {Theorem 2}Let $M$ be a connected,
holomorphically nondegenerate, algebraic
hypersurface in $\Bbb C^{ N}$. Then there exists a positive
integer $\ell$ with $1 \le \ell \le N-1$ such that if 
$H$ is a  CR mapping of class
$C^\ell$ from $M$ to another algebraic hypersurface
$M'$ in $\Bbb C^{ N}$  with  {\rm Jac } $H
\not\equiv 0$ on any open subset of $M$, then 
$H$ is locally algebraic.  Moreover, if the Levi form of $M$ is
nondegenerate at some point, then one can take $\ell=1$. 
\endproclaim

In fact in Section 1 we define a new invariant $\ell$ for any
connected real analytic hypersurface $M$ which satisfies the
conditions of Theorem 2 if $M$ is algebraic. 

Since a connected real analytic hypersurface in $\Bbb C^{2 }$ is
holomorphically nondegenerate if and only if it is not Levi flat,
the following is an immediate corollary of Theorems 1 and 2.  

\proclaim {Corollary 2} Let $M$ and $M'$ be two algebraic
hypersurfaces in $\Bbb C^{ 2}$ and assume that $M$ is connected
and not Levi flat. 
\item {\rm (i)} If $H$ is a smooth CR mapping
from $M$ to $M'$ with  {\rm Jac }$H \not\equiv 0$,  then $H$
extends holomorphically in an open neighborhood of $M$ in $\Bbb
C^{2}$.
\item {\rm (ii)} If $H$ is a CR mapping
from $M$ to $M'$ of class $C^1$ with  {\rm Jac }$H \not\equiv
0$ on any open subset of $M$,  then 
$H$ is locally algebraic.
\endproclaim 

The following is a refinement of Theorem 1 in which $H$ is assumed
to have only a previously prescribed number of derivatives,
depending only on $M$ and $M'$. The {\it degree} of an algebraic
hypersurface is the total degree of the irreducible real polynomial
defining $M$.

\proclaim {Theorem 3}Let $M$ and $M'$ be two  algebraic
hypersurfaces in $\Bbb C^{ N}$ of degrees $d$ and $d'$
respectively, and assume that
$M$ is connected and holomorphically nondegenerate. Then there
exists a positive integer $k=k(d,d',N)$ (depending only on $d,d'$
and
$N$) such that if 
$H$ is a  CR mapping from $M$ to $M'$ of class $C^k$ with  {\rm Jac
}$H
\not\equiv 0$, then 
$H$ extends holomorphically in an open neighborhood of $M$ in
$\Bbb C^{ N}$ and is algebraic.
\endproclaim

The following shows that the condition of holomorphic
nondegeneracy is necessary for the holomorphic extendability of
CR mappings to hold.

\proclaim {Theorem 4}Let $M$ be a connected real analytic
hypersurface in $\Bbb C^{N }$  which is holomorphically degenerate
at some point $p_1$.  Let
$\po\in M$ and suppose  there exists a germ at $\po$ of a smooth
CR function on $M$ which does not extend holomorphically to any
full neighborhood of $\po$ in $\Bbb C^{N }$.  Then there exists a
germ at $\po$ of a smooth CR diffeomorphism from $M$ into itself,
fixing $\po$, which does not extend holomorphically to any
neighborhood of
$\po$ in $\Bbb C^{N }$.  
\endproclaim

Our final result deals with analytic extendability of CR
mappings between hypersurfaces in complex spaces of different
dimensions. Let $M$ be a real analytic hypersurface, $p \in M$, and
$\rho$ a defining function of $M$ in a neighborhood of $p$. Recall
[D1], [D2], [Le] that if
$M$ does not contain a complex
analytic variety of positive dimension through point
$p$ then there exists $C> 0$ such that for any complex analytic
curve parametrized by $Z=\gamma(t)$ with
$\gamma(0)=p$, 
$$\od(\rho(\gamma(t),\overline{\gamma(t)}))\le C\ 
\od(\gamma(t)), \tag 0.1 $$ where
$\od(\rho(\gamma(t),\overline{\gamma(t)}))$ and $\od(\gamma(t))$
denote the orders of vanishing of
$\rho(\gamma(t),\overline{\gamma(t)})$ and $\gamma(t)$,
respectively, at
$t=0$. In this case we let $m_p$ be the
smallest integer for which (0.1) is satisfied with $C=m_p$.
\proclaim {Theorem 5}
 Let $M\subset \Bbb C^{N }$ be an
algebraic hypersurface. Assume that there is no nontrivial complex
analytic variety contained in $M$ through $\po \in M$, and let
$m=m_{\po}$ be the integer defined as above. If
$H:M \to S^{2N+1}\subset {\Bbb C}^{N+1 }
$ is a CR map of class $C^{m}$, where $S^{2N+1 }$ denotes the
boundary of the unit ball in ${\Bbb C}^{N+1 }$, then
$H$ admits a holomorphic extension in a neighborhood of $\po$.
\endproclaim

The paper is organized as follows.  In Section 1 we introduce a
new invariant for real analytic hypersurfaces, which will be used
in the proofs of Theorems 1 and 3 and could also be of
independent interest.  The proofs of Theorems 1 and 2 are given in
Sections 2 and 3 respectively. Section 4 is devoted to the proof
of Theorem 3.  In Section 5, we study properties of families of
CR automorphisms for holomorphically degenerate hypersurfaces and
give a proof of Theorem 4.  The proof of Theorem 5 is
given in Sections 6 and 7.

The authors are grateful to Leonard Lipshitz for suggesting to them
the proof of Lemma 4.2.

\heading  \bf \S 1 A new invariant for real analytic
hypersurfaces 
\endheading 

Let $M$ be a real analytic hypersurface in $\Bbb C^{N }$ through
$0$ and
$\po
\in M$ close to $0$.  If
$\rho(Z,\zzb)$ is a defining function for $M$ near $0$, with
$\rho(\po, \bar \po) = 0$ and $d\rho (\po, \bar \po)\not= 0$, we
define the {\it Segre surface} through $\bar \po$ by
$$ \Sigma_{\po }= \{ \z \in \Bbb C^{ N}: \rho(\po, \z) = 0 \}.
$$ 
Note that $\spo$ is a germ of a smooth holomorphic hypersurface in
$\cn$ through $\bar \po$. Let $L_1, \ldots, L_n$, $n = N-1$, 
given by $L_j = \sum_{k=1}^Na_{jk}(Z,\Z)
\prt {\Z_k} $ be a basis of the CR vector fields on $M$ near
$0$ with the $a_{jk}$ real analytic.   If
$X_1,
\ldots, X_n$,
 are the complex vector fields given by $X_j =
\sum_{k=1}^Na_{jk}(\po,\z)
\prt {\z_k} $, $j =1,\ldots, n$, then $X_j$ is tangent to $\spo$
and the $X_j$ span the tangent space to
$\spo$ for $\z$ in a neighborhood of $\bar \po$, with $(\po,\z)
\mapsto a_{jk}(\po,\z)$ holomorphic near $0$ in $\Bbb C^{2N }$.    For a
multi-index
$\alpha = (\a_1,\ldots,\a_n)$ we define
$c_\a(Z,\po,\z)$ in $\Bbb C\{Z,\po,\z\}$, the ring of convergent
power series in $3N$ complex variables, by 
$$
c_\a(Z,\po,\z)= X^{\a} \rho(Z+\po, \z), \tag 1.1
$$
where $X^\a = X_1^{\a_1}\cdots X_n^{\a_n}$.

Note that since the $X_j$ are tangent to $\spo$, we have
$c_\a(0,\po,\z) = 0$  for all $\po \in M$ and $\z \in \spo$ in a
neighborhood of $0$. In particular, $c_\a(0,\po,\bar\po) = 0$
for $\po \in M$ close to $0$.   We say that $M$ is {\it essentially
finite} at
$\po$ if the ideal
$(c_\a(Z,\po,\bar
\po))$ generated by the $c_\a(Z,\po,\bar
\po)$, $\a \in \Bbb Z_+^n$, in the ring $\Bbb C\{Z\}$ is of finite
codimension.  (By the Nullstellensatz, this is equivalent to the
condition that the functions $Z \mapsto c_\a(Z, \po, \bar \po)$
have only $0$ as a common zero near the origin for $\po$ fixed
and $\a \in \Bbb Z_+^n$.) This definition of essential
finiteness, which does not depend on either the choice of
holomorphic coordinates or that of the defining function, coincides
with that given in [BJT] and that given in [BR2] in a slightly
different form. The present definition has the advantage of
avoiding the use of the implicit function theorem, thus making 
explicit calculations easier.   

We introduce here a new invariant which will give us a bound on the number of
derivatives needed in Theorem 2. If $M$ is essentially
finite at
$\po \in M$ fixed as above, let $ \ell(\po)$ be the
minimum positive integer for which the ideal generated by
$\{c_\a(Z, \po, \bar \po): |\a| \le \ell(\po)\}$ is of
finite codimension in $\Bbb C \{Z\}$. It follows from the
definition of essential finiteness and the fact that $\Bbb
C\{Z\}$ is a Noetherian ring that $\ell(\po)$ is finite. 
It is clear that $\ell(\po) \ge 1$.  

\proclaim {Proposition 1.2} Let $M$ be a connected real analytic
hypersurface which is holomorphically nondegenerate. Then there
is an integer $\ell(M)$, with $ 1\le \ell(M)\le N-1$ such that
$\ell(p) =\ell(M)$ for all
$p$ in an open, dense subset of $M$.  Moreover, $\ell(M) =1$ if
and only if $M$ is generically Levi nondegenerate.    
\endproclaim
\demo {Proof}
We need to introduce the following vector-valued functions. For a
multi-index $\a$, Let $V_\a$ be the real analytic function
defined near $0$ in
$\Bbb C^{N }$ by 
$$ V_\a (Z,\Z)= L^\a \rho_Z(Z,\Z), \tag 1.3
$$
where $\rho_Z$ denotes the gradient of $\rho$ with respect to
$Z$.  In the rest of the proof we shall say that a property
holds {\it generically} on $M$ (or an open neighborhood  of $\po$
in $M$) if it holds in an open, dense subset of $M$. 

\proclaim {Lemma 1.4}
 Let $M$ be a connected real analytic 
hypersurface in $\Bbb C^{ N}$.  Then $M$ is  holomorphically
nondegenerate if and only if $\{ V_\a(Z,\Z), \a \in
\Bbb Z^{ N}_+\}$ span $\Bbb C^{ N}$ generically  in a
neighborhood of $0$ in $M$. 
 \endproclaim

\demo {Proof}  We note first that the condition that the $V_\a$
span
$\Bbb C^{ N}$ is independent of the choice of coordinates and
defining function. We
introduce here normal coordinates near $0$, $Z=(z,w)$,  $z\in \Bbb C^{
n}$,
$w\in \Bbb C^{ } $, such that $M$ is given there by 
$$ w = Q(z,\-z, \-w),\ \  \hbox {with}\   \ Q(z,0,w) \equiv
w,$$     (or equivalently  by $\wb =
\Qb(\zb, z, w)$).  We put $\po = (\zo,\wo)$; we may then take
$$
X_j = \prt {\chi_j} + \Qb_{\chi_j}(\chi,\zo,\wo)\prt {\t} ,  \tag
1.5
$$
where $\z=(\chi,\t)$, $\chi= (\chi_1,\ldots,\chi_n)$. Hence
$c_\a(Z,\po,\z)$ of (1.1) is given by
$$
c_\a(Z,\po,\pob) = -\Qb_{\ch^\a}(\zob, \zo,\wo) +
\Qb_{\ch^\a}(\zob, z+\zo,w+\wo). \tag 1.6
$$ 
Similarly, we have by using (1.3),
$$
V_\a(\po,\pob) =-\Qb_{\zb^\a,Z}(\zob,\zo,\wo). \tag 1.7
$$
 Hence $c_{\a,Z}(0,\po,\pob) =
-V_\a(\po,\pob)$.  If the rank of the $V_\a(Z,\zzb)$ is less than
$N$ generically, then at any point of  maximal rank $\po$ near $0$
in
$M$, by the implicit function theorem, there is a complex curve
$Z(t)$ such that $c_\a (Z(t), \po, \pob) = 0$ for all small $t$
and all $\a$.  Hence $M$ is not essentially finite at $\po$. 
Since  the set of essentially finite points is open [BR2] (and
the set of points of maximal rank is open and dense), there
exists an essentially finite point if and only if the generic
rank of the $V_\a(Z,\zzb)$ is $N$.  By [BR6], the existence of an
essentially finite point is equivalent to holomorphic
nondegeneracy of $M$.  This completes the proof of the lemma.
$\square$
\enddemo
We shall also need the following lemma, whose simple proof is
left to the reader.
\proclaim {Lemma 1.8} Let $f$ be a holomorphic function defined
in an open set $\Omega$ in $\Bbb C^{ p}$, valued in $\Bbb C^{N
}$.  If the  $\partial ^\a f(\ch)$, $\a \in \Bbb Z^{p
}_+$ span $\Bbb C^{N }$ generically in $\Omega$, then the $\partial
 ^\a f(\ch)$, $|\a| \le N-1$ also span $\Bbb C^{N }$ generically
in $\Omega$.   
\endproclaim
\proclaim {Lemma 1.9} Let $M$ be a holomorphically nondegenerate
real analytic hypersurface of $\Bbb C^{ N}$ through $0$.  There
exists an integer $k$, with $1 \le k\le N-1$ so that
$\{V_\a(Z,\Z), |\a| \le k\}$ span $\Bbb C^{N }$ generically in a
neighborhood of $0$ in $M$.  
\endproclaim

\demo {Proof of Lemma 1.9}  By Lemma 1.4, we may find $\po =
(\zo,\wo)
\in M$ so that the vectors $V_\a(\po,\pob)$ span $\Bbb C^{ N}$ as
$\a$ varies over all multi-indices.  We put
$f(\ch) =
\Qb_Z(\ch,\zo,\wo)$. Now by (1.7) and Lemma
1.8, we conclude that there exists, $1 \le k \le N-1$ such that 
the vector-valued functions
$\Qb_{{\ch^\a},Z}(\ch,\zo,\wo)$, $|\a| \le k$ span $\Bbb C^{N }$
generically for $\ch$ in a neighborhood of $0$ in $\Bbb C^{N-1
}$. This is equivalent to the nonvanishing of an $N\times N$ 
determinant $\Delta (\ch, \po) $.  We claim that the functions
$\Qb_{{\ch^\a},Z}(\zb,z,w)$,
$|\a| \le k$ also span $\Bbb C^{N }$ generically for $(z,w) \in
M$ near $0$.  For this, it suffices to show that $\Delta (\zb, Z)$
does not vanish identically for $Z \in M$ near $0$. Indeed, if  
$\Delta (\zb, Z) \equiv 0$ on $M$, then by complexifying the
variables, we would also have $\Delta (\ch, Z) \equiv 0$ for
$\ch$ near $0$ in $\Bbb C^{N-1 }$ and $Z$ near $0$ in $\Bbb
C^{N }$, contradicting $\Delta (\ch, \po) \not\equiv 0$.
$\square$
\enddemo
We return now to the proof of Proposition 1.2. We first note
that the function $\ell(p)$ is upper
semi-continuous on $M$, i.e. $\ell(p) \le \ell(\po)$ for $p$
near $\po$.  By (1.6), (1.7) and the implicit function
theorem, if follows that if $\{V_\a(\po,\pob), |\a| \le k\}$
span $\Bbb C^{ N}$, then $\ell(\po) \le k$.  Conversely, if $k$
is the smallest integer for which $\{V_\a(\po,\pob), |\a| \le k\}$
span $\Bbb C^{ N}$  generically for $\po$ in a neighborhood of
$0$, then it cannot happen that $\ell(\po) <k$ for any $\po$
near $0$.  For if so, by going to an arbitrarily close point
$p$ where the rank of
$\{V_\a(p,\bar p), |\a| \le \ell(\po) \}$ is maximal and applying
the implicit function theorem, we would obtain a complex curve of
common zeroes for the functions $\{c_\a(Z,p,\bar p): |\a| \le
\ell(\po)|\}$. This would be a contradiction, since, by the above,
$\ell(p)\le
\ell(\po)$.  This proves that the minimum of $\ell(p)$ in a
neighborhood of $0$ is the same as the smallest integer  $k$
satisfying the conclusion of Lemma 1.9. Since $M$ is connected
and real analytic, It suffices to take
$\ell(M) $ to be the smallest integer $k$ of Lemma 1.9.

To complete the proof of Proposition 1.2, it remains to show
that $\ell(M) = 1$ is equivalent to $M$ being generically Levi
nondegenerate.  For this note that an easy row and column
manipulation show that  $$\hbox{\rm det}\ [\rho_Z (Z,
\Z),L_1\rho_Z(Z,\Z),
\ldots,L_{N-1}\rho_Z(Z,\Z)]$$ is a nonvanishing multiple of the
usual Levi determinant of $M$ at $Z$.  (See e.g. [W1].)   
$\square$
\enddemo

 Proposition 2.1 and its proof suggest the
definition of a new invariant which refines the notion of
holomorphic nondegeneracy.  Let $M$ be a real analytic
hypersurface in $\Bbb C^{N }$  and $\rho$ a local defining
function.  We say that $M$ is {\it $k$-holomorphically
nondegenerate} at $Z \in M$ if the $L^\a\rho_Z(Z,\Z) $,
$|\a| \le k$, span $\Bbb C^{ N}$  with $k$ minimal. It follows
from the proof of Proposition 1.2 that if $M$ is holomorphically
nondegenerate then generically
$1\le k\le N-1$. Also, 
$M$ is $1$-holomorphically nondegenerate at
$Z$ if and only if the Levi form of $M$ is nondegenerate at $Z$.
Note that if $M$ is connected and holomorphically nondegenerate
then there exists $\ell = \ell(M)$,  $1 \le \ell(M)\le N-1$, such
that
$M$ is $\ell$-holomorphically nondegenerate at every point in
any open dense subset of $M$.  This number $\ell(M)$ is given by
Proposition 1.2; we shall call  $\ell(M)$ the {\it Levi type} of
$M$.

\heading  \bf \S 2 Proof of Theorem 1
\endheading 
First recall  from  [BR6, Theorem 2] that since $M$
is connected and holomorphically nondegenerate, the set of points
at which $M$ is essentially finite (see [BJT] and [BR1] for
definition) is not empty. On the other hand, if $M$ is
essentially finite at $p$, then $M$ is also of  finite type (in
the sense of Bloom-Graham [BG]) at $p$.  Let 
$$ U = \{p\in M: M \ \hbox {is of finite type at  } \ p\}. \tag
2.1
$$
Since $M\backslash U$ is a real analytic subset of $M$, it follows
from the above that $M\backslash U$ is a proper (possibly empty)
real analytic subset of
$M$, and hence $U$ is an open, dense subset of $M$. More
precisely, $\partial U =  M\backslash U$ is a smooth complex
hypersurface in $\Bbb C^{ N}$.  Indeed, this can be seen by
 using a theorem of Nagano [N]; if
$M\backslash U$ is nonempty, then it is given locally as the real
analytic manifold whose complexified tangent space is spanned by
the CR tangent vectors and their complex conjugates. Note that if
$\partial U$ is nonempty then it is of codimension $1$ in
$M$.   

\proclaim {Proposition 2.2} Let $M$, $M'$, and $H$ be as in Theorem
1.  Then {\rm Jac} $(H)$ does not vanish
identically on any open set in $M$.
\endproclaim

We shall need the following in the proof of the proposition.

\proclaim {Lemma 2.3}  Let $M$ be a real analytic
hypersurface in
$\Bbb C^{N}$, and assume that 
$U$ defined by {\rm (2.1)} is nonempty.  If $f$ is a continuous CR
function on $M$ and $f$ vanishes identically in some neighborhood
of
$\po$ in $U$, then $f$ vanishes identically in the whole connected
component of $\po$ in $U$.
\endproclaim 

\demo {Proof} Let $U_0$ be the connected component of $\po$ in
$U$,  and let
$$S = \{Z \in U_0: f|_V \not\equiv 0 \ \hbox {for
any neighborhood }V  \hbox { of} \ Z \ \hbox {in } \ U_0\}.   
 $$
We claim that $S$ is open and closed in $U_0$.  Indeed, it is
immediate from the definition that $S$ is closed.  To show that
$S$ is open, we let $q \in S$ and choose a connected neighborhood
$W
\subset U_0$ of $q$ sufficiently small such that $f$ extends to
one side of $M$ with boundary $W$.  (The extendability of CR
functions at
$q$ to one side of $M$ follows from the fact that $M$ is of
finite type at $q$ [BT].) If $f$ were to vanish on an open
subset of $W$, then the holomorphic function extending $f$ would
vanish identically, and hence $f$ would vanish on $W$,
contradicting the assumption that $q $ is in $S$.  This shows $W
\subset S$ and completes the proof of Lemma 2.3.
  $\square$ 
\enddemo

In the following we shall write $J(Z)$ for Jac $H(Z)$ for $Z \in
M$.

\proclaim {Lemma 2.4 } Under the assumptions of Theorem 1, if
$J|_{U_0}
\not\equiv 0$, where $U_0$ is a connected component of $U$, then
$J|_{U_0}$ is algebraic.   
\endproclaim 

\demo {Proof} Since $M$ is holomorphically nondegenerate, as
noted before by Theorem 2 of [BR6] the set where 
$M$ is essentially finite is nonempty.  Hence by Proposition 1.12
of [BR2], the set of essentially finite points in $M$ is open
and dense.  By the continuity of $J$ we may find  $p_0 \in
U_0$ such that $M$ is essentially finite at $p_0$ and
$J(p_0)\not=0$.  By a result in  [BJT] we conclude that $H$
extends holomorphically to an open neighborhood $ \co$ of $p_0$ in
$\Bbb C^{N }$. Denote by $\Cal H$ this holomorphic extension.   We
may now use Theorem 1 in [BR6] to conclude that  $\Cal H$ is
algebraic in $\co$. Since the derivatives and products of
algebraic functions are again algebraic, the holomorphic
extension $\Cal J$ of $J$ to $\co$ is also algebraic in
$\co$. Let $P(Z,X)$ be the polynomial, with polynomial
coefficients, such that $P(Z, \Cal J(Z))\equiv 0$ in $\co$. 
On the other hand, since $J$ is a CR function on $M$, we
conclude that $f(Z) = P(Z,  J(Z))$, $Z \in M$, is also CR on $M$. 
By Lemma 2.3, $f(Z) $ vanishes identically on $U_0$, since it
vanishes on $\co \cap M$.  This shows that $J|_{U_0}$ is
algebraic.
 $\square$ 
\enddemo

\proclaim {Lemma 2.5} Let $g$ be a smooth CR function on $M$ and
assume that  $g|_{U_0}$ is algebraic, where
$U
_0$ is a connected component of $U$ (given by {\rm (2.1)}).  If
$g|_{U_0}\not\equiv 0$, then $g$ cannot vanish to infinite order
at any point in the closure of $U_0$.
\endproclaim

\demo {Proof} Let $P(Z,X)$ be a polynomial such that
$$P(Z,g(Z))|_{U_0} \equiv 0.
\tag 2.6$$ 
 Suppose $P(Z,X) =
P_1(Z,X)P_2(Z,X)$, where $P_j(Z,X)$ are polynomials of positive
degree in $X$. Then since $P_j(Z,g(Z))$ is CR on $M$, $j=1,2$, we
may use Lemma 2.3 to conclude that either $P_1(Z,g(Z))|_{U_0}
\equiv 0$ or
$P_2(Z,g(Z))|_{U_0} \equiv 0$.  Hence we may assume that the
polynomial $P(Z,X) = \sum_0^k a_j(Z)X^j$ in (2.6) is
irreducible, and, in particular, $a_0(Z) \not\equiv 0$. If
$g(Z)$ vanishes of infinite order in the closure of $U_0$, it
would follow from (2.6) that the restriction of $a_0(Z)$ to
$M$ also vanishes of infinite order at that point.  Since
$a_0(Z)$ is a polynomial, and $M$ is real analytic, this
would imply $a_0(Z)$ vanishes identically, contradicting the
irreducibility of $P$.    
  $\square$ 
\enddemo

\demo {Proof of Proposition {\rm 2.2}}  Define $E $ by
$$E = \{Z \in M: J|_V \not\equiv 0 \ \hbox {for
any neighborhood }V  \hbox { of} \ Z \ \hbox {in } \ M\}.   
 $$
Since $E$ is nonempty by assumption, the proposition will follow
from the connectedness of $M$ if we prove that $E$ is open and
closed.  The closedness of $E$ is immediate from the definition. 
We shall show that $E$ is open.   First, if $\po \in E\cap U$,
and  $U_0$ is the connected component of
$U$ containing $\po$, then by
Lemma 2.1 we have $U_0 \subset E$ . If $\po \in E\cap \partial U$, 
 and $V$ is a sufficiently small neighborhood of
$\po$, $V$ will intersect at most two connected components of $U$,
say
$U_1$ and $U_2$.  (For this, recall that $\partial U$ is a
smooth submanifold of $M$ of codimension $1$.)  By definition of
$E$, either $J|_{U_1\cap V} \not\equiv 0$ or $J|_{U_2\cap V}
\not\equiv 0$.  It then follows from Lemmas 2.4 and 2.5 that $J$
cannot vanish to infinite order at $\po$; therefore $J|_{U_j\cap
V} \not\equiv 0$, $j=1,2$.  By Lemma 2.1, this shows $V \subset
E$, which completes the proof that $E$ is open and the proposition
follows. $\square$ 
\enddemo

We shall need the following result, which is probably known. 
(See also Lemma 6.1 of [BJT] for a special case of this result.) 

\proclaim {Lemma 2.7 } Let $G(z,w)  $ be a holomorphic function
in a neighborhood of $0$ in $\Bbb C^{p+1 }$ with $G_w(z,w)
\not\equiv 0$. Let
$f$ be a smooth  function defined in a neighborhood of $0$ in $\Bbb
R^{ p}$ satisfying
$$  G(x,f(x)) \equiv 0,  \tag 2.8
$$ 
for $x$ in a neighborhood  of $0$ in $\Bbb R^{p }$.  Then
$f$ is real analytic in a neighborhood of $0$. 
\endproclaim

\demo {Proof} We wish to apply the following theorem due to
Malgrange [M]: Let $Y$ be the germ of a real analytic set in
$\Bbb R^{q }$  through 0 containing a germ $\Sigma$ of a smooth
manifold through $0$.  If
$\hbox {dim}_\Bbb R   Y = \hbox {dim}_\Bbb R  \Sigma$, then
$\Sigma$ is real analytic.

We write $G(x,s) = \sum_{j=0}^{\infty} a_j(x)s^j$, where each
$a_j$ is a convergent power series.    We may assume that the
$a_j(x)$,
$j = 0,1, \ldots,$ have no common factors (as convergent power
series in
$p$ variables). We define the real analytic set  $Y \subset \Bbb
R^{ p+2} $ as follows.  Let $Q_1(x,y), Q_2(x,y)$ be the real valued
functions determined by
$$ \sum_0^\infty a_j(x)(y_1+iy_2)^j = Q_1(x,y)+iQ_2(x,y)
$$
with $y = (y_1,y_2) \in \Bbb R^{2 }$, and let $Y$ be the germ of
the real analytic set at $0$ defined by
$$Y = \{(x,y)\in \Bbb R^{p+2 }:Q_1(x,y)=Q_2(x,y)=0\}.$$    Let
$\Sigma$ be the germ at $0$ of the smooth submanifold of $Y$ given
by the parametrization
$$ \Sigma = \{ (x,g_1(x),g_2(x)): x \in \Bbb R^{p } , \ \hbox
{where } \ f(x) =g_1(x)+ig_2(x) \}.$$
Clearly $\hbox {dim}_\Bbb R  \Sigma = p$; the desired real
analyticity of $f$ will follow from Malgrange's result above if
$\hbox {dim}_\Bbb R  Y = p$.  To prove this last equality, note
that if $(x,y) \in Y$, then each $ y_j$ is determined by $x$ up to
finitely many values unless all the $a_j$ vanish at $x$.  However,
since the $a_j(x)$ have no common factors, we claim that the
dimension of their common zeros is less than or equal to $p-2$. 
For this, note first by the Noetherian theorem, there exists $k$
such that
$\{ x: a_j(x) = 0, j=0,1,\ldots,\} = \{ x: a_j(x) = 0,
j=0,1,\ldots,k\}$. Now the claim can be seen by expressing
the
$a_j$ as   Weierstrass polynomials with respect to the same variable and
applying an elimination method  as e.g. in Lemma 5.1 of [BJT].   
  $\square$ 
\enddemo

\demo {Proof of Theorem 1} Let $U$ be given by (2.1).  Since
by Proposition 2.2 $J$, the Jacobian of $H$, does not vanish
identically on any open subset of $M$, we conclude as in the
proof of Lemma 2.4 that $H$ is algebraic on each connected
component of $U$. In order to show that $H$ extends
holomorphically to a neighborhood of $M$, by standard arguments it
suffices to show that $H$ is real analytic in a neighborhood of
each point in
$M$. Let $p_0 \in M$.  We claim that $H$ is algebraic in some
neighborhood of $\po$.  Indeed, if $\po \in U$, then by the
above, $H$ is algebraic in the component of $\po$.  If $\po \in
\partial U$, then $\po$ is in the boundary of at most two
components, say $U_1$ and $U_2$.  Hence for $j = 1,\ldots, N$ there
exist polynomials $p_{1j}(Z,X)$ and $p_{2j}(Z,X)$ such that
$$p_{kj}(Z,H_j(Z))|_{U_k} \equiv 0,\ \  k = 1,2, \ \ j =
1,\ldots,N. \tag 2.9
$$ Let $p_j(Z,X) = p_{1j}(Z,X) p_{2j}(Z,X)$. It follows from (2.9)
that $p_j(Z,H_j(Z))|_{U_1 \cup U_2} \equiv 0$, which proves the
claim. By taking a real analytic parametrization of $M$, we may
apply Lemma 2.7 to conclude that $H$ is real analytic at every
point.  This completes the proof of Theorem 1.
$\square$ 
\enddemo

The following example shows that the assumption that $M$ is algebraic
cannot be dropped.

\remark {Example {\rm 2.10}} The following example was given by
Ebenfelt in [E]. Let $t=\theta(\xi,s)$ be the unique solution of the
algebraic equation $\xi(t^2+s^2) - t = 0$, with       
$\theta(0,0)=0$.  Let $M$ and $M'$ be the  hypersurfaces through $0$ in
$\Bbb C^2$ given respectively by
$$\Im w= \theta(\arctan |z|^2, \Re w)  \ \ \ \ \ \ \ \ \  
\Im w'= (\Re w') |z'|^2. $$
(Note that $M$ and $M'$ are both of infinite type at $0$.) Let
$H=(f,g)$ be the mapping given by $f(z,w) = z$, $g(z,w) =
e^{-(1/w)}$ for $\Re w >0$ and $g(z,w) \equiv 0$ for $\Re w \le
0$.  It is shown in [E] that $H$ is a smooth CR mapping defined
in a neighborhood of $0$ which maps $M$ to $M'$.  However, it is
clear that $H$ does not extend holomorphically to a neighborhood
of $0$ in $\Bbb C^{ 2}$.  Note that $M'$ is algebraic, but $M$
is not.
\endremark 

\heading  \bf \S 3 Proof of Theorem 2 
\endheading

To prove Theorem 2, we take $\ell=\ell(M)$ to be the Levi type of
$M$ defined at the end of Section 1.  Let
$U$ be the open set in $M$ defined in (2.1) and
$\partial U$ its boundary.  Since Jac $H$ does not vanish
identically on any open set, it does not vanish identically on any
connected component of $U$.  We fix such a component $U_0$ and choose
$\po \in U_0$ such that 
Jac $H(\po) \not= 0$,
and  $\ell(\po) = \ell(M)$; the latter is possible by Proposition
1.2. Note that in particular $M$  is essentially finite at $\po$. 
As before, we assume that local normal coordinates on $M$ and $M'$
are chosen so that
$\po= 0$ and $H(\po) = 0$.  In these coordinates we write
$H=(f,g), f=(f_1,\ldots, f_n)$, $n = N-1$.

Using the methods of proof of Lemma 6.1 of [BR1] and Proposition
2.5 of [BR6], we obtain the following.  For each
$j$, $1
\le j
\le n$, there exists a positive integer $N_j$ and  algebraic
functions 
$a_{jk}(u^\gamma_p,v^\beta)$, $|\gamma|, |\beta| \le \ell$, $0 \le k
\le N_{j-1}$,
$1\le p
\le n$,    holomorphic
near
$u^{\gamma} _{p,0} = L^{\gamma}  \fb_p(0)$ and $v^\beta_0 = 0$,
such that we have in a neighborhood of $0$ in $M$:
$$ f_j^{N_j} + \sum_{k=0}^{N_j-1}a_{jk}(
L^{\gamma} \fb_p,L^\beta \gb)f_j^k\equiv 0, \ \ \ \ j=1,
\dots,n.
\tag 3.1
$$ 
(Although in [BR1] and [BR6] the mapping $H$ was assumed to be
smooth, the proof of (3.1) uses derivatives only up to length
$\ell$; hence (3.1) holds also when $H$ is only of class $C^\ell$.)
We may now move to a point $p_1$, arbitrarily close to $0$ near which
the roots of the polynomials (3.1) are analytic functions of the
coefficents.  We conclude that near that point $p_1$
$$
f_j = \Psi_j(
L^{\gamma} \fb_p,L^\beta \gb),\ j = 1,\ldots, n,  \tag 3.2
$$
with $\Psi_j$ analytic.   By the standard use of the reflection
principle, as for instance in [BR1] (see also Proposition 7.1 below),
we conclude that the
$f_j$ extend holomorphically near $p_1$.  It then follows easily that
the same holds for $g$ near $p_1$. We continue to denote by $H = (H_1,
\ldots, H_N)$ the original CR map as well as its holomorphic
extension in a neighborhood of $p_1$. We may now apply Theorem 1 of
[BR6] to conclude that $H$ is algebraic in a neighborhood in $\Bbb
C^{ N}$. That is, for $j = 1,\ldots, N$, there exist polynomials
$P_j(Z,X)$ with holomorphic polynomial coefficients such that
$P_j(Z,H_j(Z)) \equiv 0$ is a neighborhood of $p_1$.  Let $k_j(Z) =
P_j(Z,H_j(Z))|_M$; each $k_j$ is a CR function on $M$ which
vanishes on a neighborhood of $p_1$ in $M$.  Hence, by the argument
of the proof of Theorem 1, $k_j$ must vanish identically in the
connected component $U_0$ of $p_1$ in $U$.  This shows that the
restriction of $H_j$ to each connected component of $U$ is
algebraic.  It remains to show the same holds near a point $\po \in
\partial U$.  Since $\partial U$ is a smooth hypersurface of $M$,
$\po$ is in the closure of two connected components of $U$. For each
$j$ we take the product of the two polynomials corresponding to the
two connected components to obtain an algebraic equation satisfied
by $H_j$ in the closure of the union of these components.  This
completes the proof of Theorem 2. $\square$

\remark {Example {\rm 3.3}} In Theorem 1 it was sufficient to
assume that Jac $H$ does not vanish identically on $M$.  The following
example shows that for the conclusion of Theorem 2 to hold, unlike
that of Theorem 1, we must assume the stronger condition that Jac 
$H\not\equiv 0$ on each component of $U$.  In the smooth case (i.e.
Theorem 1), we show in the course of the proof that the two
conditions are actually equivalent.  
For this, let $M$ and $M'$ be the  hypersurfaces through $0$ in
$\Bbb C^2$ given respectively by
$$\Im w= (\Re w)^3 |z|^2,   \ \ \ \ \ \ \ \ \  
\Im w'= 2\frac {(\Re w') \theta ( |z'|^2)} { 1-[\theta (
|z'|^2)]^2},$$
where $\theta (t)$ is the unique real solution vanishing for $t=0$,
of the polynomial equation $X^3+X -t=0$.  Note that both $M$ and
$M'$ are algebraic and generically Levi nondegenerate. Consider the
CR mapping
$H$ defined on $M$   by  $H=(f, g)$ with 
$$
\left\{\eqalign{f(z,w)&=zw, \ \ \ \Re w \ge 0\cr f(z,w)&=zw\exp w,
\ \ \ \Re w \le 0 \cr}
\right\},  \qquad \left\{\eqalign{g(z,w)&=w^2, \ \ \ \Re w \ge 0\cr
g(z,w)&=0,
\ \ \ \Re w \le 0 \cr}
\right\}.
$$
The reader can easily check that $H$ is of class $C^1$ on $M$ and
that
$H$ maps  a
neighborhood of $0$ in $M$ into $M'$.  Note that $H$ is  not algebraic, and
Jac $H(z,w)\not=0$ for $(z,w) \in M,
\hbox{with}\ 
\Re w > 0$, but Jac $H(z,w) \equiv 0$ for $(z,w) \in M,
\hbox{with}\ 
\Re w < 0$.
Note that in this example the number $\ell$ given by Theorem 2 is $1$.
 \endremark 
\medskip

\heading  \bf \S 4 Proof of Theorem 3
\endheading

In this section we shall indicate the modifications to the proof of
Theorem 1 needed to  give Theorem 3. For a polynomial 
$$P(Z,X)= a_J(Z)X^J +a_{J-1}(Z)X^{J-1} +\ldots a_0(Z) \tag 4.1
$$ with
polynomial coefficents $a_j(Z)$, where               $a_J(Z)
\not\equiv 0$, by the {\it total degree} of
$P$ we shall mean the total degree of $P$ as a polynomial in the
variables
$(X,Z)$. If $f(Z)$ is an algebraic function, by the {\it total degree}
of $f$ we shall mean the minimum of the total degrees among
polynomials $P(Z,X)$ for which $P(Z,f(Z))\equiv 0$.  We need the
following lemma, whose proof is based on the Artin Approximation
Theorem and was suggested to us by Leonard Lipshitz.

\proclaim {Lemma 4.2 }For any positive integers $d$ and $n$ there
exists a positive integer $k = k(d,n)$ such that if $f$ is a function
of class $C^k$ defined in a neighborhood of $0$ in $\Bbb R^{ n}$
and satisfies a nontrivial polynomial equation $p(x,f(x)) \equiv 0$,
where
$p(x,Y)$ is a polynomial of $n+1$ variables of total degree $\le
d$, then $f$ is real analytic in a neighborhood of $0$.
\endproclaim
\demo {Proof} Let $\Delta
(x)$ be the discriminant of $p(x,Y)$ regarded  as a polynomial in the
indeterminate $Y$. By eliminating repeated factors in the
factorization of
$p(x,Y)$, we may assume that $\Delta (x) \not\equiv 0$. We use the
following consequence of the Artin Approximation Theorem [A1],
[A2], [BDLV]: 

{\sl Let $p(x,Y)$, $x \in \Bbb R^{n }$,  be a polynomial in $Y$
with polynomial coefficients.  Then for any positive integer $r$ there
exists a positive integer $k$ (depending only on $r$, $n$ and the
total degree of $p$) such that if
$f_1(x)$ is a formal series for which
$p(x,f_1(x)) = O(|x|^k)$ there is a convergent series $g(x)$ such
that $p(x,g(x)) \equiv  0$ and $g(x) - f_1(x) = O(|x|^r)$. }

In fact, the statement above is a special case of Theorem 6.1 of
[A2] or Theorem 3.2 of [BDLV]. We shall apply the above with
$r=(d_1 + d_2(J(J-1))+ 1)/2$, where
$d_1$ is the degree of the  polynomial
$\Delta(x)$ defined above, $d_2$ is the degree of $a_J(x)$.  Note
that $r$ is bounded by an expression  depending only on the total
degree of
$p$.  If $k$ is  given by the statement above, then we claim that the
conclusion of Lemma 4.2 holds with this choice of $k$. For this,
let
$f_1(x)$ be the truncated Taylor series of $f(x)$ up to degree $k$
and 
 let
$g(x)$ be the convergent series given by the statement above.  
If $\rho(x)$ is a root of $p(x,Y)$, then $a_J(x)\rho(x)$ is a root
of the monic polynomial 
$$q(x,Y) =  Y^J+ a_{J-1}(x)Y^{J-1}+
\ldots + a_J^{J-1}(x)a_0(x). \tag 4.3 $$
 In particular, $a_J(x)f(x)$ and
$a_J(x)g(x)$ are both roots of $q(x,Y)$. If $g(x) \not\equiv f(x)$,
let $\tau_1(x)= a_J(x)f(x)$, $\tau_2(x)= a_J(x)g(x)$, and let 
$\tau_3(x),\ldots
\tau_J(x)$ be the rest of the roots of
$q(x,Y)$ (counted with multiplicity). Then the discriminant of
$q(x,Y)$ is $$a_J^{J(J-1)}(x)\Delta(x) = 
-(a_J(x)f(x)-a_J(x)g(x))^2\Pi
(\tau_j(x)-\tau_k(x)), \tag 4.4$$ where the indices on the right
hand side run over $j \not= k$ and either $j$ or $k$ is not equal to
$1$ or
$2$. Since $q(x,Y)$ is a monic polynomial, the $\tau_k$ are bounded.
Hence the right hand side of (4.4) vanishes to order at
least $2r$.  On the other hand, since the left hand side of (4.4) is
of degree
$\le d_1 +d_2(J(J-1))$, both sides must vanish identically, by the
choice of
$r$, contradicting the assumption that $\Delta(x) \not\equiv 0$. 
This contradiction shows that $g(x) \equiv f(x)$, which completes
the proof of the lemma.
  $\square$ 
\enddemo

We shall now prove Theorem 3.  We start with the following analogue of 
Proposition 2.2. 

\proclaim {Lemma 4.5} If $k$ is sufficiently large, and Jac $H
\not\equiv 0$, then Jac
$H$ does not vanish identically on any open set in $M$.
\endproclaim  
\demo {Proof}
 By the
connectedness of $M$, and using Lemma 2.3, it suffices to show that
if Jac $H \not\equiv 0$ in some connected component $U_0$ of $U$,
then Jac $H \not\equiv 0$ on any component of $U$ which is
contiguous to $U_0$.  As in \S 2 we may find $p \in U_0$ at which
$M$ is essentially finite and  Jac $H(p) \not= 0$. Hence, the
components of $H$ extend holomorphically in a neighborhood $V$ of $p$
in $\Bbb C^{N }$ . By [BR6], there are polynomials $P_j(Z,X)$ of
the form (4.1) such that the $j$th component of $H$ satisfies
$P_j(Z,H_j(Z))
\equiv 0$ for $Z \in V$. By the proof given in [BR6], one can see
that the total degree of the $P_j(Z,X)$ is bounded by a number
which depends only on the total degrees of the defining functions
of $M$ and $M'$.  Hence Jac $H$ is the root of a polynomial
$P(Z,X)$ whose total degree is bounded by a number depending only
on the total degrees of the defining functions of $M$ and $M'$. 

 By propagation of the
zeroes of CR functions in $U_0$, 
$$P(Z, \hbox {Jac} \ H(Z)) \equiv 0 \tag 4.6$$
for $Z \in U_0$, and, as in \S 2, we may assume that $P(Z,X)$ is
irreducible.  If Jac
$H$ were to vanish on a component
$U_1$ contiguous to
$U_0$, then it vanishes to order at least $k-1$ on the boundary
between $U_0$ and $U_1$.  Hence the constant coefficient $a_0(Z)$ of $P$
must also vanish to order $k-1$ there. Since the degree of $a_0(Z)$
is bounded by the total degree of $P$, we would have $a_0(Z) \equiv 0$
if $k-1$ is greater than the total degree of $P$.  Since this would
contradict the irreducibility of $P$, we conclude as before that Jac
$H$ does not vanish on any open set in $M$.
$\square$ 
\enddemo

We shall now show that each component $H_j$ of $H$ is algebraic with
total degree bounded by a number depending only on the total  degrees
of the defining functions of $M$ and $M'$.  By Lemma 4.2, this will
complete the proof of Theorem 3. As in the argument above, for a
fixed component $U_0$ of $U$, there exists polynomials
$P_j(Z,X)$, with total degree bounded by  number depending only on the
total degrees of the defining functions of $M$ and $M'$, such that
$P_j(Z,H_j(Z))\equiv 0$, $j = 1,\ldots,
N$ for $Z \in U_0$.  By the connectedness of
$M$, it suffices to show that $P_j(Z,H_j(Z))\equiv 0$,  in any
component $U_1$ of $U$ adjacent to $U_0$. By Lemma 4.5, Jac $H
\not\equiv 0$ on $U_1$, so that one can find polynomials $\tilde
P_j(Z,X)$, with total degree bounded by a number depending only on the
total degrees of the defining functions of $M$ and $M'$, such that    
 $\tilde P_j(Z,H_j(Z))\equiv 0$, $j = 1,\ldots,
N$ for $Z \in U_1$. We now have that
$$  P_j(Z,H_j(Z))\tilde P_j(Z,H_j(Z))\equiv 0,\  j = 1, \ldots, N, \tag
4.7$$
for $Z \in U_0\cup U_1$. Hence if $k$ is sufficiently large
(depending only on the total degrees of $M$ and $M'$, we may apply
Lemma 4.2 to conclude that $H$ is real analytic in the interior of
the closure of $U_0\cup U_1$. By unique continuation of analytic
functions, it follows that $P_j(Z,H_j(Z)) \equiv 0$ for $Z \in U_0\cup
U_1$. This completes the proof of Theorem
3. $\square$
\medskip
It should be noted that in general the integer $k$ in Theorem 3
could be much larger that $\ell(M)$, the Levi type of $M$ defined at
the end of Section 1. The following example shows that in Theorem 3,
the integer
$k$ cannot be taken to be
$\ell(M)$. 

\remark {Example {\rm 4.8}}  Let $M$ and $M'$ be given as in Example
3.3.  Consider the CR mapping $H$ defined on $M$ by 
$H=(f, g)$ with $f(z,w) = zw$ and $g(z,w) = w^2$ for $\Re w \ge
   0$, and $g(z,w) = -w^2$ for $\Re w \le 0$.  As observed by Peter
Ebenfelt, $H$ is of class
$C^1$ and  Jac $H$ does not vanish identically on any open subset of
$M$. However, $H$ is algebraic but does not extend holomorphically in
any neighborhood of $0$ in $\Bbb C^{2 }$ Note that $\ell(M)=1$
here.   
\endremark    

\heading  \bf \S 5 Families of CR automorphisms; Proof of Theorem
4
\endheading

We shall prove Theorem 4 in this section. By Proposition 4.2 of
[BR6], since $M$ is holomorphically degenerate at $p_1$, it is
holomorphically degenerate at $\po$ also. We choose local
coordinates near $\po$ for which $\po$ is the origin, and let
$X= \sum_1^N a_j(Z)\prt
{Z_j}
$, where the $a_j(Z)$ are germs at $0$ of holomorphic functions,
be a nontrivial holomorphic vector field tangent to $M$ near $0$. 
Let $h(Z)$ be a smooth CR function defined on $M$ near $0$ which
does not extend holomorphically to any neighborhood of $0$ in $\Bbb
C^{N }$. We may choose $h$ so that $h(0) = 0$.   Let $\phi(t,Z)$
be the flow of 
$X$ for
$t
\in
\Bbb C^{ }
$,
$|t|$ small, i.e. $\phi(t,Z)$ satisfies the holomorphic ordinary
differential equation
$$
\prt t  \phi(t,Z) = A(\phi(t,Z)), \  \
\phi(0,Z) = Z, \tag 5.1
$$
where $A= (a_1,\ldots,a_{N})$. Let $Y$ be the complex vector
field on $M$ near $0$ obtained from $X$ by multiplication of the
coefficients by
$h$ i.e.,
$$ Y = hX= \sum_1^N h(Z)a_j(Z)\prt
{Z_j}    . \tag 5.2
$$
\proclaim {Lemma 5.3} The differential equation
$$\prt t K(t,Z) = h(\phi(K(t,Z),Z)), \ \ \ \ K(0,Z) \equiv 0, \tag
5.4
$$
has a smooth solution $K(t,Z)$ defined for $(t,Z)$ in a
neighborhood of
$(0,0)$ in
$\Bbb C^{ } \times M$ with $t \mapsto K(t,Z)$ holomorphic in
$t$ for fixed $Z$, and $Z\mapsto K(t,Z)$  CR for $t$ fixed.  
\endproclaim 
\demo {Proof} Let $F(t,Z) = h(\phi(t,Z))$.  Note that since $h(0)
= 0$, it follows that $F(0,0) = 0$.  Since $\phi(t,Z)$ is
holomorphic in $Z$ and $h$ is CR, $Z \mapsto F(t,Z)$ is also CR. 
We claim that $t \mapsto F(t,Z)$ is holomorphic.  Indeed,
let $\tilde h$ be any smooth extension of $h$ to a neighborhood
of $0$ in $\Bbb C^{N }$.  By the chain rule,
$$ \prt {\bar t} (\tilde h(\phi(t,Z),Z)) = \tilde h_Z \cdot \prt
{\bar t}
  \phi(t,Z) + \tilde h_{\bar Z}  \cdot \overline {\prt { t}
\phi(t,Z)}  .\tag 5.6
$$ 
The first term on the right in (5.6) is zero since $\phi$ is
holomorphic in $t$. By (5.1), the second term is $\tilde h_{\bar
Z}
\cdot
\bar A$, which equals $(\bar X h)(\phi(t,Z),Z)$, since $\tilde h =
h$ on
$M$ and $\bar X$ is tangent to $M$. This term also vanishes since
$\bar X$ is a CR vector field and $h$ is a CR function. The
existence of a smooth solution $K(t,Z)$, holomorphic in $t$ for
$Z \in M$ is then given by the holomorphic theory of ordinary
differential equations.

To see that $K(t,Z)$ is CR, let $L$ be a CR vector field near
$0$.  Then using (5.4), a simple calculation shows that $LK(t,Z)$
satisfies the ODE 
$$\prt { t} LK(t,Z) = L \prt { t}   K(t,Z) =
(Xh)(\phi(K(t,Z),Z))LK(t,Z),
\tag 5.7
$$
with $LK(0,Z) \equiv 0$.  By uniqueness, $LK(t,Z) \equiv 0$ for all
$t$. 
$\square$ 
\enddemo 
\proclaim {Lemma 5.8} Let $\psi(t,Z) = \phi(K(t,Z),Z)$, with
$K(t,Z)$ given by Lemma {\rm 5.3}, be defined for
$(t,Z)$ in a neighborhood of $(0,0)$ in
$\Bbb C^{ } \times M$.  Then $\psi$ defines a complex CR flow for
the vector field $Y$ defined by {\rm (5.2)}.  That is,
$$
\prt { t}   \psi(t,Z ) = h(\psi(t,Z ))A(\psi(t,Z )) \tag 5.9
$$
with $\psi(0,Z) = Z$ and $t \mapsto \psi(t,Z )$ is holomorphic
for fixed $Z$.  Furthermore,
$Z
\mapsto
\psi(t,Z)$ is CR for each $t$.
\endproclaim
\demo {Proof} Since $\prt { t}   \psi(t,Z )= (\prt {
t}       \phi)(K(t,Z),Z)\cdot\prt {
t}       K(t,Z)$,  (5.9) and the holomorphy of $\psi(t,Z
)$ with respect to $t$ are immediate from the properties of
$K(t,Z)$ given in Lemma 5.3. The fact that $\psi(t,Z)$ is CR for
fixed $t$ follows easily since the same is true of $K(t,Z)$.
  $\square$ 
\enddemo

\proclaim {Lemma 5.10} Let $R(t,Z)$ be a smooth function defined
in $\de \times V$ with  $\de =
\{ t \in \Bbb C^{ } : |t| < \epsilon\}$ and $V$ a neighborhood of
$0$ in $M$.  Assume that $R$ is holomorphic in $t$ for fixed $Z$,
$R(t,0) \equiv 0$, and for each $t \in \de$ there exists
$\Cal O_t$, a neighborhood of $0$ in $\Bbb C^{ N}$ such that
$Z \mapsto R(t,Z)$ extends holomorphically to $\Cal
O_t$.  Then there exist $\eta > 0$, $t_0 \in \de$, and $\Cal O$ an
open neighborhood of $0$ in $\Bbb C^{N }$,  such that
$R(t,Z)$ extends holomorphically to $\{t: |t-t_0|< \eta\} \times
{\Cal O}$.  
\endproclaim
\demo {Proof} For a positive integer q let $E_q \subset \de$ be
given by 
$$E_q = \{ t
\in
\de: | {D^\alpha} R(t,Z)| \le \alpha !q^{|\alpha|}\ \hbox {for
all}\ Z \in M,\ |Z| < {1\over q}, \ \hbox {and all} \ \alpha \}, 
$$
where $D^\alpha$ denotes differentiation on $M$ in some fixed local
parametrization of $M$ near $0$.   Since by assumption $\cup_q E_q =
\de$, and the
$E_p$ are closed, we may apply the Baire Category Theorem to find
$q_0$ such that $E_{q_0}$ has nonempty interior. That is, there
exist $t_0 \in \de$, $\eta >0$, and $C>0$ such that for all nonzero
$\alpha$
$$
|{D^\alpha} R(t,Z)| \le \alpha ! C^{|\alpha|} \ \hbox {for
all}\ Z \in M,\ |Z| < {1\over C},\ |t-t_0|<\eta.
$$
It follows that $R(t,Z)$ extends continuously to a neighborhood of
the form
$\{t:|t-t_0|<
\eta\}
\times {\Cal O}$, separately holomorphic in $Z$ and $t$. The
lemma is then a consequence  of Hartog's Theorem. 
  $\square$ 
\enddemo 

\demo {Proof of Theorem 4} We shall prove the theorem by
contradiction. Suppose that for every germ of a smooth self CR map
of $M$ fixing
$0$ there exists a neighborhood of $0$ in $\Bbb C^{N }$ to
which the map extends holomorphically.  In particular, this would
imply that for each $t\in \Bbb C^{ } $ small, the CR map $Z
\mapsto \psi(t,Z)$ extends holomorphically to a neighborhood of
$0$, where $\psi$ is given by Lemma 5.8. Since $\psi(t,Z) =
\phi(K(t,Z),Z)$, we may apply Lemma 2.7 to conclude that for each
$t \in \Bbb C^{ }$ sufficiently small, $Z \mapsto K(t,Z)$ extends
holomorphically to a  neighborhood of $0$ in $\Bbb C^{N }$.  We
may now apply Lemma 5.10, to conclude that $K(t,Z)$ extends holomorphically to $\{|t-t_0|< \eta\} \times
{\Cal O}$. In particular, we conclude, by differentiating in $t$
and using (5.4),  the function
$Z \mapsto h(\phi(K(t_0,Z),Z))$ extends holomorphically near the
origin.  Since $K(t, 0) \equiv 0$ (by uniqueness in (5.4)) and
the map $Z \mapsto \phi(K(t_0,Z),Z)$ is a local biholomorphism
near the origin, we conclude that $h(Z)$ extends holomorphically
in a neighborhood  of $0$ in $\Bbb C^{N }$.  This
contradicts the assumption on $h$ and completes the proof of 
Theorem 4.
  $\square$ 
\enddemo

\heading  \bf \S 6 Properties of mappings into the sphere
\endheading

We will prove Theorem 5 in this section and the next.  We begin
with some notation.  We write $H = (H_1, \ldots, H_{N+1})$, and, as
before, $N=\no$.  After a local holomorphic change of variables, we may
assume that
$S^{2N+1}$ is given by the defining function 
$$Z'_{\nt} +\bar
{Z}'_\nt + \sum_{j=1}^{\no} |Z'_j|^2.   \tag 6.1$$

\proclaim {Lemma 6.2} Under the assumptions of Theorem 5, $M$ is
pseudoconvex in a neighborhood  of $\po$ and there exist points of
strict pseudoconvexity arbitrarily close to $\po$.
\endproclaim 
\demo {Proof} It follows from the hypothesis that $M$ is of
Bloom-Graham finite type near $\po$; hence $H$ extends
holomorphically to one side of $M$ near $\po$ [BT], [T]. Let  
$$\rho^*(Z,\Z)=H_{\nt}(Z)+\overline{H_\nt(Z)}+\sum_{j=1}^{n+1}
|H_j(
Z)|^2.\tag 6.3$$
Hence $\rho^*$ is defined on one side of $M$, of class $C^m$ up to
$M$, and is plurisubharmonic.
 For any small analytic
disc $A(\z)$ attached to $M$ near $\po$, the function
$\z \mapsto \rho^*(A(\z), \overline {A(\z)})$ is subharmonic in the
unit disc and vanishes on its boundary.  By the maximum
principle, the interior of the disc maps to the unit ball in $\Bbb
C^{N+1 }$.  Since the discs cover one side of $M$, we
can apply  the Hopf Lemma to conclude
$\prt {\nu }  {\rho^*(Z,\Z)} |_{\po}
\not= 0$, where
$\nu$ is the normal direction to $M$ at $\po$.
Hence $\rho^*$ is a plurisubharmonic
 defining function for $M$, proving the pseudoconvexity of $M$.

Suppose that there is no strongly pseudoconvex point near $p_0$. Then,
by the semi-continuity of the counting function of the positive
eigenvalues of the Levi form, we can assume that for some point
$p_1$ close to
$\po$ the number of the non-zero eigenvalues of the Levi form of $M$
there attains a local maximum value $r < n$. Now, by using a local
change of coordinates near $p_1$, we may assume that $p_1=0$ and 
$M$ can be given near this point by the following equation
$$Z_{\no}+\Z_{\no}=\sum_{j=1}^{r}|Z_j|^2+h(Z_1,\cdots,
Z_{\no}),$$  where $h(Z)=O(|Z|^3)$.
Consider the hyperface $M^*\subset \Bbb C^{\no-r } $ defined by
$Z_\no+\Z_\no =h(0,\cdots,0,Z_{r+1},\cdots, Z_\no)$. Since it
does not contain any non-trivial analytic variety inside (for,
otherwise, $M$ cannot be of D'Angelo finite type), one sees that it
cannot be Levi-flat and thus  its Levi form has a positive eigenvalue
at a point $w^*$ near $0$ in $M^*$. Then  
the Levi form of $M$ has at least $r+1$ 
positive eigenvalues at $(0,w^*)$. This contradicts the definition
of $r$ and completes the proof of Lemma 6.2.
  $\square$ 
\enddemo

An immediate consequence of Lemma 6.2 with Theorem 2 in [H2] is that
$H$ is an algebraic map. This fact will be useful later.

As in Section 1, we choose normal coordinates, $Z=(z,w)$, for $M$ 
vanishing at $\po$ so that $M$ is given by an equation of the form:
$$t=\phi(z,\bar z,s) \tag 6.4$$
where $w=s+it$ 
and $\phi(0,\bar z,s)\equiv\phi(z,0,s)\equiv 0$.
Let
$L_1,\ldots, L_n$ be a basis for the CR vector fields on $M$ near
$0$ given by
$$L_j = \prt {\zb_j}   -2i{ \phi_{\zb_j} \over 1 +i\phi_s} \prt 
{\wb}  , \  \ j = 1, \ldots, n.  \tag 6.5
$$
 Note that if  
$L^\a = L_1^{\a_1}\ldots L_n^{\a_n}$, then 
$$  L^\a|_0 = \prtf {^{|\a|}} {\zb^\a} |_0.  \tag 6.6
$$

We will also assume $H(0)=0$.

\proclaim {Lemma 6.7} After a  rotation of the vector
$(H_1,\cdots,H_\no)$, one can find a sequence of $n$ multi-indices 
$(\beta^1,\cdots,\beta^{n})$ with $|\beta^j|\le m$
such that for $j = 1,\ldots, n$,
$${ L}^{\b^j}\H_j|_0\not=0,\ \ \h{\rm but}\ {
L}^{\beta^j}\H_\ell |_0= 0\ \ \h{\rm for}\ 1\le j< \ell \le
n. 
\tag 6.8$$ 
\endproclaim 
\demo {Proof} 
We note first that by applying $L^\b$ to (6.3) we have
$$ L^\b \H_\nt (0) = 0,\ \h {\rm for all} \ |\b| \le m.\tag 6.9 
$$
 We first show
that there exists a multi-index
$\b$, $|\b| \le m $, and an integer $k \le \no$ such that
$L^{\b}\H_k(0)
\not=0$.  We argue by contradiction.  If no such $\b$ exists, then
by (6.6) we have $\H_j(Z) = O(|w|) + o(|z|^m)$, $j=1,\ldots,\no$.
However, since (6.3) is a defining function for
$M$, $t-\phi(z,\bar z,s) = \rho^*(Z,\Z)h(Z,\Z)$ for some
nonvanishing $h$.  Using (6.9), we have 
$
H_\nt = aw + o(|w| + |z|^m), 
$ 
with $a\not= 0$. Combining these gives
$$ \rho^*(Z,\Z) = aw +\bar a\wb + o(|w| + |z|^m).  \tag 6.10
$$
Then the complex analytic variety $V^* = \{Z = (\z,\ldots,\z,0):
\z \in \Delta\}$ has order of contact at least $m+1$  with $M$ at
$0$, contradicting the definition of $m$ given by (0.1).

Without loss of generality we may assume that there exists
$\b^1$, $|\b^1| \le m$ such that $L^{\b^1}\H_1(0)\not= 0$ and
$L^{\b^1}\H_j(0)= 0$, $ 1 < j \le \no$.
Next, we show that there exists a multiple index $\beta^2$ with 
$|\beta^2|\le m$ so that ${ L}^{\beta^2}\H_j(0)\not=0$ for
some $1 < j \le \no$. Indeed, if this is not the case, then
$H_j=O(|w|)+o(|z|^m)$, $j = 2,\ldots, \no$. Write
$H_1(Z)=P_1(z)+O(|w|)+o(|z|^m)$ with $P_1(z)$ a polynomial
in $z$ of degree $\le m$ and $P_1(0)=0$. Then using again (6.9)
we obtain 
$$\rho^*(Z,\Z)=aw+\bar{a}\wb+|P_1(z)|^2+o(|w|+|z|^m).
$$
Let $V$ be the complex analytic variety defined by $w=0$ and
$P_1(z)=0$. Without loss of generality, we may assume 
that
$P_1(b_1\tau,\cdots,b_{n-1}\tau,z_{n}) =\sum _{j=0}^N
a_j(\tau)z_{n}^j$ with $a_j(\tau)\not\equiv 0$ for some $j>0$ and
an (n-1)-tuple $(b_1,\cdots,b_{n-1})$.
 By the Puiseux expansion, there exists an $N^*>>1$
such that $z_{n}(\tau^{N^*})$ is a holomorphic function in $\tau$,
$z_{n}(0)=0$, and
$P_1(b_1\tau^{N^*},\cdots,b_{n-1}\tau^{N^*},z_{n}(\tau^{N^*}))\equiv
0$. Thus, we obtain a holomorphic curve $\gamma(\tau)$ defined by
$$z_1=b_1\tau^{N^*},\cdots,z_{n-1}=b_{n-1}\tau^{N^*},z_{n}=z_{n}(\tau^{N^*}),
w=0.$$
If $N'= \od(\gamma(\tau))$, the order of vanishing of
$\gamma(\tau)$ at 0, then since
$\rho^*(\gamma(\tau),\overline{\gamma(\tau)})=o(|\gamma(\tau)|^m)$
we have
$$\od(\rho^*(\gamma(\tau),\overline{\gamma(\tau)})
 >m \ \od (\gamma(\tau)),$$
which again would contradict the definition of $m$. Hence $\b^2$
must exist. Now, applying a suitable rotation to
$(H_2,\cdots,H_n)$, we may assume that
${L}^{\beta^2}\-H_2|_0\not=0$ for some $|\beta^2|\le m$ but
${L}^{\beta^2}\-H_l|0=0$ for each $l>2$.
Arguing inductively, we obtain the  proof of Lemma 6.7.     
  $\square$ 
\enddemo

We shall use Lemma 6.7 to obtain equations for the components
$H_i$ of the mapping $H$.
\proclaim {Proposition 6.11} Let $\b^j$ and $H_i$ be as in Lemma
{\rm 6.7}, and let
$$V(Z,\Z)=\pmatrix{L}^{\beta^1}\H_1&{L}^{\beta^1}\H_2&
\ldots&{L}^{\beta^{1}}\H_{n}\cr\vdots&\vdots&\vdots\cr
{L}^{\beta^{n}}\H_1&
{L}^{\beta^{n}}\H_2&\ldots
&{L}^{\beta^{n}}\H_{n}\cr
\endpmatrix.$$
Then  $V(0,0)$ is invertible, and if
$$\xi=V^{-1}\pmatrix{L}^{\beta^1}\H_\nt\cr\vdots\cr
{L}^{\beta^{n}}\H_{n+2}\cr\endpmatrix,\ 
\eta=V^{-1}\pmatrix {L}^{\beta^1}\H_{\no}\cr\vdots\cr
{L}^{\beta^{n}}\H_{\no}\cr\endpmatrix, \
F=\pmatrix H_1\cr\vdots\cr H_{n}\cr \endpmatrix, \tag 6.12$$
then the following holds on $M$:
$$F= -{\xi}-{H_\no}{\eta}.
\tag 6.13$$ 
\endproclaim

\demo {Proof} Apply  $L^{\b^j}$, $j = 1, \ldots, n$, to (6.3),
and then solve the resulting system of linear equations  for $
F$.  
  $\square$ 
\enddemo

\heading  \bf \S 7 Proof of Theorem 5
\endheading

We shall complete the proof of Theorem 5 in this section.   The main
step will be to prove that $H$ is meromorphic. Then the result in
Chiappari [Ch] will give the desired holomorphic extension. It will be
convenient to have the following criterion.
\proclaim{ Proposition 7.1} Let $M$ be a minimal algebraic
hypersurface in
$\Bbb C^{N } $ ($N>1$) defined near $0$,
 and let $k(Z,\Z)$ be a vector-valued algebraic continuous CR function
defined on $M$ near $0$. Assume that $h(Z,\Z)$ is also a continuous
 algebraic CR function on $M$ near $0$ such that
$Q_2(Z,\-Z,\-{k(Z, \Z)})h(Z, \Z)=Q_1(Z,\-Z,\-{k(Z,\Z)})$
 for $Z \in M$ near $0$, where
$Q_j$ ($j=1,2$) are holomorphic algebraic functions near
$(0,0,\-{k(0,0)})$ such that $Q_2(Z,\-{Z},\-{k(Z,\Z)})\not \equiv 0$
near the origin in $M$. Then
$h(Z,\Z)$ admits a meromorphic extension near $0$ in $\Bbb C^{N }$. 
Moreover, when
$Q_2\equiv 1$, then $h$ admits a holomorphic extension near $0$.
\endproclaim
\demo {Proof} We will use the edge of wedge theorem for the proof.
Assume that $M$ is given in normal coordinates $Z=(z,w)$ by  equation
(6.4) and that each CR function
defined near $0\in M$ can be extended to the side $D^+$ given by $t >
\phi(z,\-z,s)$. For each nonzero vector  $v \in \Bbb C^{n } $, let
$M_v=\{(xv,s+i\phi(xv,x\-{v}, s)): x\in \Bbb R^{n } ,\ s\in \Bbb R^{
} \}$. Since
$Q_2(Z,\-{Z},\-{k(Z,\Z)})$ cannot vanish identically in an open set
of $M$
 near the origin, we can assume for some fixed $v_0$,
$Q_2(Z,\-{Z},
\-{k(Z,\Z)})\not \equiv 0$ on $M_{v_0}$ in any neighborhood of $0$.
Now we define
$G:\ \Bbb C^{\no } \rightarrow \Bbb C^{\no }$, a local
biholomorphism with
$G(0) = 0$ by
$$G(z,w)=(zv_0,w+\phi(zv_0,z\-{v_0},w§)).$$ 
Note that $G(\Bbb R^{\no }) \subset M_{v_0} $
near $0$. Furthermore, $G$  maps a standard straight wedge $W^+$
with edge
$\Bbb R^{\no } $ near $0$ to $D^+$. Write $W^-=\-{W^+}$, the complex
conjugate of
$W^+$.
Denoting the holomorphic extensions of $h$ and $k$ to $D^+$ by
the same letters, let
 $$\Phi^+(Z)=h\circ G(Z), \ \ \ Z\in W^+,$$ 
$$\Phi^-(Z)={Q_1(G(Z),\-{G(\-Z)},\-{k(G(\-Z),\-{G(Z)})})\over
Q_2(G(Z),\-{G(\-Z)},\-{k(G(\-Z),\-{G(Z)})})},  \ \ \ Z\in
W^-.$$  Then, by our hypothesis, one sees that
$\Phi^-$ is a meromorphic algebraic function in $W^-$. Thus,
 $\sum_{j=1}^{N}a_j(Z)\Phi^-(Z)^j\equiv 0$., where $a_j$'s are
polynomials. Notice that $a_N\Phi^-$ satisfies a polynomial equation
with leading coefficient 1. We conclude that $a_N\Phi^-$ is bounded in
$W^-$. By the Riemann removable singularity theorem, it
holds that
$a_N\Phi^-$ extends holomorphically to $W^-$. 

Away from a
proper real analytic subset of $\Bbb R^{ \no} $, we have
$a_N(X)\Phi^+(X)=a_N(X)\Phi^-(X)$ for $X$ near $0$ in  $\Bbb R^{\no
}$, by our construction. From the classical edge of the wedge
theorem,  it follows that $a_N(Z)\Phi(Z)$ can be extended
holomorphically to an open subset of $0$. Thus
$\Phi^+$ has a meromorphic extension near $0$. Now, since $G$ is
a local  biholomorphic map, we conclude that $h$ extends
meromorphically across
$0$. 

Finally, if $Q_2\equiv 1$, then $\Phi^-(Z)$ is  bounded
near
$0$ for $Z\in W^-$. The above argument shows that $h$ is holomorphic
near $0$ in this case. This completes the proof of the proposition.
 $\square$ 
\enddemo

We now begin the proof of Theorem 5.  We proceed according to the
following two cases:
\roster
\item
"{Case I}": ${L}_j(\eta(Z,\Z))\equiv 0$, $j=1,\ldots, n$ for $Z$ in
$M$ near
$0$, i.e. the components of the vector $\eta$ are all CR functions.
\item "{Case II}":
At least one of the components of $\eta(Z,\-Z)$ is not a CR
function in any neighborhood of $O$ in $M$.
\endroster

We first assume the hypothesis of Case I and prove that $H$ extends
meromorphically 
 near
$0$.  Since $F$ and $H_\no$ are CR functions, we conclude from
(6.13) that $\xi$ is also a CR function.  Applying the
last part of Proposition  7.1 to the first two equations of
(6.12), we conclude that $\xi$ and $\eta$  both extend
holomorphically to a full neighborhood of $0$ in $\Bbb C^{\no }$ and
hence are both real analytic near $0$ in $M$.  

Rewriting (6.3) we obtain on $M$,
$$ H_\nt +\-H_\nt +F\cdot\-F + H_\no\-H_\no = 0. \tag 7.2$$
Replacing $F$ in (7.2) by using (6.13), we have 
$$ H_\nt +\-H_\nt +H_\no a +\-H_\no\-a + |H_\no|^2 + c = 0, \tag
7.3$$
where 
$$ a = \eta \cdot \-\eta, \ \  b = 1 + |\eta|^2, \ \ c = |\xi|^2 \tag
7.4$$
Rewriting (7.3) we have
$$ H_\nt +\-H_\nt +H_\no (a +\-H_\no b) + \-H_\no\-a + c = 0, \tag
7.5$$
Applying $L_j$, $j=1,\ldots, n$ to (7.5) we obtain
$$ L_j\-H_\nt +H_\no L_j(a +\-H_\no b) + L_j(\-H_\no\-a) + L_jc = 0,
\tag 7.6$$
We consider now two cases.
\roster
\item
"{Case Ia}": $L_j(a +\-H_\no b)\equiv 0$, $j = 1,\ldots, n$
\item "{Case Ib}":
For some $j$,  $L_j(a +\-H_\no b)\not\equiv 0$. 
\endroster

In Case Ia, we conclude that $H_\no = -a/b$ and hence
$H_\no$ extends holomorphically by Proposition 7.1.  (Note that $b$ is
nowhere vanishing.)  Hence from (6.13) and then (7.2) it follows
that $H_j$, $j=1,\ldots,\nt$ also extend holomorphically.

In Case Ib we apply Proposition 7.1 to equation (7.6) to conclude
that $H_\no$ extends meromorphically, and hence again by (6.13) and
(7.2) we conclude that all the $H_j$ extend meromorphically.

We now consider Case II.  Then choose  $j, \ell \in
\{1,2,\ldots n\}$ so that
$L_j
\eta_\ell \not\equiv 0$ in any neighborhood of $0$ in $M$. Applying
$L_j$ to the $\ell$th component of (6.13), we obtain
$$L_j\xi_\ell + H_\no L_j \eta_\ell = 0 \tag 7.7$$  We then use
Proposition 7.1 to conclude that $H_\no$ has a meromorphic
extension.  Making use first of (6.13) and then of (7.5), we
conclude that all the components of $H$ extend meromorphically.

We have now proved that $H$ extends to a meromorphic mapping in a
neighborhood of $0$.  Since $H$ maps one side of $M$ to the ball
and maps $M$ to the sphere, we may now use Theorem 1 of Chiappari
[Ch], (generalizing the result of Cima-Suffridge [CS]) to conclude
that $H$ extends holomorphically.  This completes the proof of
Theorem 5.

\enddocument